\documentclass{amsart}
\usepackage{amssymb,amscd}

\newtheorem{lemma}{Lemma}
\newtheorem{prop}[lemma]{Proposition}
\newtheorem{cor}[lemma]{Corollary}
\newtheorem{thm}[lemma]{Theorem}

\newtheorem{thm?}[lemma]{Theorem?}
\newtheorem{ques}[lemma]{Question}
\newtheorem{conj}[lemma]{Conjecture}
\newtheorem{prob}[lemma]{Problem}
\newtheorem*{claim}{Claim}

\DeclareFontEncoding{OT2}{}{} 
  \newcommand{\textcyr}[1]{%
    {\fontencoding{OT2}\fontfamily{wncyr}\fontseries{m}\fontshape{n}%
     \selectfont #1}}
\newcommand{\Sha}{{\mbox{\textcyr{Sh}}}}

\title{Abelian points on algebraic curves}
\author{Pete L. Clark}
\address{The Mathematical Sciences Research Institute \\
17 Gauss Way \\ Berkeley, CA 94709} \email{plclark@msri.org}

\newcommand{\F}{\ensuremath{\mathbb F}}
\newcommand{\Fp}{\ensuremath{\F_p}}

\newcommand{\Fq}{\ensuremath{\F_q}}

\newcommand{\Q}{\ensuremath{\mathbb Q}}
\newcommand{\R}{\ensuremath{\mathbb R}}
\newcommand{\Z}{\ensuremath{\mathbb Z}}

\newcommand{\ra}{\ensuremath{\rightarrow}}

\newcommand{\ord}{\operatorname{ord}}

\newcommand{\PP}{\ensuremath{\mathbb P}}

\newcommand{\Gal}{\operatorname{Gal}}

\newcommand{\Br}{\operatorname{Br}}

\newcommand{\ab}{\operatorname{ab}}

\newcommand{\Qp}{\mathbb{Q}_p}

\renewcommand{\div}{\operatorname{div}}

\begin{document}
\maketitle

\begin{abstract}
We study the question of whether algebraic curves of a given genus $g$ defined
over a field $K$ must have points rational over the maximal abelian extension
$K^{\ab}$ of $K$.  We give: (i) an explicit family of
diagonal plane cubic curves without $\Q^{\ab}$-points, (ii) for every
number field $K$, a genus one curve $C_{/\Q}$ with no
$K^{\ab}$-points, and (iii) for every $g \geq 4$ an
algebraic curve $C_{/\Q}$ of genus $g$ with no $\Q^{\ab}$-points.  In
an appendix, we discuss varieties over $\Q((t))$, obtaining in
particular a curve of genus $3$ without $(\Q((t))^{\ab}$-points.
\end{abstract}

\noindent Convention: All varieties over a field $K$ are assumed to be
nonsingular, projective and
(as is especially important for what follows) geometrically irreducible.

\section{Introduction}

\noindent In \cite{Frey}, G. Frey demonstrated the existence of an
algebraic variety $V_{/\Q}$ with no points rational over the
maximal abelian extension $\Q^{\ab}$ of $\Q$ (or ``without abelian
points.'') His argument uses a mixture
of elliptic curve theory and valuation theory;\footnote{See also
the ``elliptic curve free'' proof given in \cite{FJ}.} from it,
one can deduce the existence of an abelian variety $A_{/\Q}$ and a
principal homogeneous space $V_{/\Q}$ for $A$ such that
$V(\Q^{\ab}) = \emptyset$.  But one does not have any bound on the
dimension of $A$, nor any information about the curves lying on
$V$.\footnote{One does, of course, know that there exist curves on $V$, but
regarding the \emph{existence} of an object as \emph{information} about it
seems to invite philosophical controversy.}
\\ \\
The purpose of this paper is to take a closer look at curves
without abelian points. First, we shall give ``optimally''
concrete and simple examples of varieties without abelian points.
For this, let us not delay the statement (and proof!) of the
following
\begin{thm}
\label{EASY} Let $p \equiv -1 \pmod 3$ be a prime number, and let
$a,\ b, \ c$ be integers which are prime to $p$.  The curve
\begin{equation}
\label{ONE} aX^3 + bpY^3 + cp^2Z^3 = 0 \end{equation} has no
$\Q_p^{\ab}$-points.
\end{thm}
\noindent Proof: First observe that $C$ has no $\Q_p$-rational
points.  For, if not, we would have a solution $(x,y,z) \in
\Z_p^3$ with $\min( \ord_p(x), \ \ord_p(y), \ \ord_p(z)) = 0$ and
this is visibly impossible: looking at the equation we see that
$p$ must divide first $x$, then $y$, then finally $z$.  Indeed,
the same argument works in any finite extension $K/\Q_p$ in which
the relative ramification degree $e(K/\Q_p)$ is prime to $3$.
\\ \indent
But now suppose that there exists a solution in the ring of
integers of $\Q(\mu_N)$ for some positive integer $N$.  Write $N =
M \cdot p^i$ with $(M,p) = 1$, and let $K$ be the completion of
$\Q(\mu_N)$ at some prime lying over $p$.  We have $e(K/\Q_p) =
\varphi(p^i) = p^{i-1}(p-1)$, which is, by our assumption on $p$,
prime to $3$.  Because the maximal abelian extensions of both $\Q$
and $\Q_p$ are those generated by all roots of unity,\footnote{The
Kronecker-Weber theorem.  In fact the elementary structure theory
of finite extensions of local fields would suffice.} the proof is
complete.
\\ \\
For ``most'' number fields $K$, we can choose $p$ such that the
curves (\ref{ONE}) fail to have $K^{\ab}$-rational points.  More
precisely:
\begin{cor}
\label{EASYC} Let $K$ be a number field whose Galois closure does
\emph{not} contain $\Q(\mu_3)$.  Then there exists a prime $p$
such that the curves (\ref{ONE}) have no $K^{\ab}$-rational
points.
\end{cor}
\noindent Proof: Thanks to our assumption on $K$, Cebotarev's
density theorem guarantees the existence of infinitely many primes
$p \equiv -1 \pmod 3$ such that $p$ splits completely in $K$.  We
then have an embedding $K^{\ab} \hookrightarrow \Q_p^{\ab}$, and
by Theorem \ref{EASY} we conclude that $aX^3 + pbY^3 + p^2cZ^2 =
0$ has no $K^{\ab}$-rational points.
\\ \\
In fact we will prove the following:
\begin{thm}
\label{NUMBERFIELD} For $K$ a number field, there is a genus
$1$ curve $C_{/\Q}$ with $C(K^{\ab}) = \emptyset$.
\end{thm}
\noindent Our second goal is to determine for which genera $g$
there exists a genus $g$ curve $C_{/\Q}$ with $C(\Q^{\ab}) =
\emptyset$. Since quadratic extensions are abelian, one obvious
sufficient condition for abelian points is for $C_{/\Q}$ to admit
a degree $2$ morphism to a curve $Y$ with $Y(\Q) \neq \emptyset$.
Taking $Y = \PP^1$ we see that there are abelian points on all
hyperelliptic curves, and in particular on all curves of genus $0$
or $2$.
%
\begin{thm}
\label{HIGHERG} For all $g \geq 4$, there exists $C_{/\Q}$ of
genus $g$ and such that $C(\Q^{\ab}) = \emptyset$.
\end{thm}
\noindent Conspicuously missing is the case of $g = 3$.
Equivalently, we wonder:
\begin{ques}
\label{GENUS3} Must a nonsingular plane quartic curve $C_{/\Q}$
have an abelian point?
\end{ques}
\noindent The question remains of interest over $\Q_p$ (and especially,
over $\Q_2$) and has prompted me to begin a more systematic
study of rational points on curves over local fields.
\\ \\
The first step of the proofs has a similar flavor to the
arguments of \cite{Frey}: namely, we begin in $\S 2$ by
constructing genus one curves $C_p$ defined over $\Q_p$ (in fact,
over arbitrary $p$-adic fields) without abelian points via Galois
cohomological methods. \\ \indent The second step is to ``pull
back'' these $C_p$'s to curves $C_{/\Q}$ with $C_{/\Q_p} \cong
C_p$.  We thus construct curves $C_{/\Q}$ without
$\Q_p^{\ab}$-abelian points and \emph{a fortiori} without
$\Q^{\ab}$-abelian points. That such pullbacks should exist (for
suitably chosen $C_p$'s) seems unsurprising -- it would be strange
indeed if among members of a certain ``class'' of varieties there
existed $\Q_p$-varieties without $\Q_p^{\ab}$-rational points but
no $\Q$-varieties without $\Q^{\ab}$-rational points -- but to
prove their existence poses some technical challenges, as global
Weil-Ch\^atelet groups $H^1(\Q,E)$ are agreeably large but
unwieldy in structure.  Moreover for our subsequent applications
we need genus one curves $C_{/\Q}$ without abelian points and with
further stringent conditions on the \emph{index}, as described in
Theorem \ref{ELL}.  The proof of Theorem \ref{ELL} revisits some
ideas of \cite{CRELLE}.  In particular, as in \cite{CRELLE}, use
is made of elliptic curves $E_{/\Q}$ with rank zero and known
(finite!) Shafarevich-Tate groups.\footnote{Whereas in
\cite{CRELLE} the $E_{/\Q}$ came courtesy of the work of
Gross-Zagier and Kolyvagin, this time we find ourselves using an
$E_{/\Q}$ coming from the (slightly earlier) work of Rubin and
Mazur.}
\\ \indent The curves of Theorem \ref{HIGHERG} are constructed as degree two
coverings of the genus one curves of Theorem \ref{ELL}, following
a suggestion of B. Poonen.
\\ \\
Obviously the present work is very far from being an authoritative
treatment of abelian points on algebraic varieties.  Many deeper
questions remain, especially concerning varieties which have
points everywhere locally and varieties ``on the other side of the
Calabi-Yau line,'' e.g. hypersurfaces defined by a form of degree
$d$ in more than $d$ variables.  These two issues are related: for
instance, it is unknown whether there exists a cubic surface
$S_{/\Q}$ without $\Q^{ab}$-rational points, but by a theorem of
Lang \cite{Lang} such an $S$ must have $\Q_p^{\ab}$-rational
points for all $p$.  Further remarks on these and other related
issues are made in $\S 4$.
\\ \\
After seeing an early draft of this note, B. Poonen commented that
it might be of interest to investigate the existence of abelian
points on varieties over arbitrary fields $F$, as is done in
\cite{Pal} for solvable points.  In an appendix, we take $F =
\Q((t))$, obtaining in particular curves of genus $3$ and
geometrically rational surfaces without abelian points.
\\ \\
Acknowledgements: We are happy to acknowledge the financial
support of the Mathematical Sciences Research Institute and the
mathematical support of Jean-Louis Colliot-Th\'el\`ene and Bjorn
Poonen.

\section{Local fields}
\noindent We identify principal homogeneous spaces
$V$ for an abelian variety $A_{/K}$ with Galois cohomology classes
$\eta \in H^1(K,A)$. It is thus natural to speak of a field
extension $L/K$ such that $V(L) \neq \emptyset$ as a
\emph{splitting field} for $V$ (or for the corresponding
cohomology class).
\\ \\
For $\eta \in H^1(K,A)$, its \emph{period} is the order of $\eta$
as an element of the torsion abelian group $H^1(K,A)$ and its
\emph{index} is the greatest common divisor of all $[L:K]$ for
$L/K$ a splitting field for $\eta$.  When $A = E$ has dimension
one, the index of $\eta$ is equal to the least positive integer
which is the degree of a $K$-rational divisor on the corresponding
genus one curve $C$. To any $K$-rational divisor of degree $n \geq
3$ there corresponds a degree $n$ embedding of $C$ into
$\PP^{n-1}$ and conversely (up to linear equivalence of divisors and
automorphisms of projective space); in particular, the curves of index $3$
are precisely the plane cubic curves without rational points.
\\ \\
In this section the ground field $K$ is a finite extension of
$\Q_p$, with valuation ring $R$ and residue field $\F_q$, $q =
p^a$. Let $r$ be the cardinality of the (finite!) group of roots
of unity in $K$, so that $q-1 \ | \ r$ and $\frac{r}{q-1}$ is a
power of $p$.
\\ \\
Our point of departure is the following result:
\begin{thm}(Lang-Tate, \cite{LT})
\label{LT} Let $A_{/K}$ be an abelian variety with good reduction,
and let $V$ be a principal homogeneous space for $A$ whose order,
$n$, is prime to $p$.  For
a finite field extension $L/K$, the following are equivalent: \\
(i) $V(L) \neq \emptyset$. \\
(ii) The ramification index $e(L/K)$ is divisible by $n$.
\end{thm}
\noindent From this we deduce the following result, a sharpening
of the example of \cite[$\S 2.3$]{WCII}.
\begin{prop}
\label{LOCAL1} Let $V$ be a principal homogeneous space of an
abelian variety $A_{/K}$ with good reduction, whose order $n$ is
prime to $p$ and does not divide $r = \#\mu(K)$.  Then $V(K^{\ab}) =
\emptyset$.
\end{prop}
\noindent Proof: Suppose on the contrary that $L/K$ is an abelian
extension with $V(L) \neq \emptyset$.  We may decompose $L/K$ into
a tower of extensions
\[K = F_0 \subset F_1 \subset F_2 \subset F_3 = L \]
such that $F_1/F$ is unramified, $F_2/F_1$ is totally ramified and
of degree prime to $p$, and $F_3/F_2$ is totally ramified of
degree a power of $p$ (e.g. \cite{CL}).  Put $d = [F_2:F_1]$.  \\
\indent We work with the corresponding Galois cohomology class
$\eta \in H^1(K,A)[n]$, our assumption being that $\eta |_{F_3} =
0$.  Since $[F_2:F_3]$ is prime to $n$, we have $\eta |_{F_2} = 0$
\cite[Prop. 11]{WCII}. By Theorem \ref{LT}, we have that $n \ | \ d$, and by
the known structure of totally tamely ramified Galois extensions
of local fields \cite{CL}, there exists a uniformizer $\pi$ of
$F_1$ such that $F_2 = F_1[X]/(X^d-\pi)$.  Since $L/K$ is abelian,
so is $F_2/F_0$, which implies that $\mu_d \subset F_0$. In other
words, $n \ | \ d \ | \ r$, a contradiction.
\\ \\
It is thus in our interest to give conditions for the existence of
classes $\eta \in H^1(K,A)$ of a given order $n$ (prime to $p$ and
not dividing $r$).  To ease the notation in the proof of the
next result, we
write $B$ for the dual abelian variety $A^{\vee}$.  (Of
course $B \cong A$ for principally polarizable abelian varieties,
and in particular for elliptic curves.)
\begin{lemma}
\label{LOCAL2} Let $\tilde{A}_{/\Fq}$ be the reduction of $A$.  If
$\tilde{A}(\Fq)$ has an element of order $n$, with $(n,p) = 1$,
then there exists $\eta \in H^1(K,B)$ of order $n$.
\end{lemma}
\noindent Proof: By a seminal theorem of Tate \cite{Tate}, the
discrete abelian group $H^1(K,B)$ and the profinite abelian group
$A(K)$ are in Pontrjagin duality.  It follows that the finite
abelian groups $H^1(K,B)[n]$ and $A(K)/nA(K)$ are in duality and
hence are isomorphic.  Reduction modulo the maximal ideal of $K$
gives an epimorphism $\mathcal{R}: A(K) \ra \tilde{A}(\Fq)$, whose
kernel $K$ is uniquely $n$-divisible.  So $\mathcal{R}$ induces an
isomorphism $A(K)/nA(K) \stackrel{\sim}{\ra}
\tilde{A}(\Fq)/n\tilde{A}(\Fq)$, and the result follows.
\noindent

\begin{thm}
There is a genus one curve $C_{/K}$ with $C(K^{\ab}) = \emptyset$.
\end{thm}
\noindent Proof: Assume for the moment the following
\begin{claim}
There exists a prime $\ell$ such that \\
$\bullet$ $\ell$ is prime to $q(q-1)$. \\
$\bullet$ There is an elliptic curve $E_{/\F_q}$ with $\ell \ | \
\# E(\F_q)$.
\end{claim}
\noindent Then: let $\tilde{E}$ be any lift of $E$ to an elliptic
curve over $K$ (e.g., choose representatives in $R$ of the
coefficients of a Weierstrass equation of $E$).  By Lemma
\ref{LOCAL2} there exists a class $\eta \in H^1(K,\tilde{E})$ of
order $\ell$, and since $\ell$ is prime to $(q-1)q$ it does not
divide $r$.  So by Proposition \ref{LOCAL1}, the corresponding
genus one curve $C_{/K}$ has $C(K^{\ab}) = \emptyset$.
\\ \indent
It remains to prove the claim.  Once we refer the reader to the
Deuring-Waterhouse classification of the integers $N$ which are $\# E(\F_q)$ 
for some elliptic
curve $E_{/\F_q}$ \cite[Theorem
4.1]{Waterhouse}, she (resp. he) may prefer to work out the proof
for herself (resp. himself).  Our proof will use the fact that there exist 
such $N$ of the form $q+1-t$ when $t$ is an integer of absolute value at 
most $2\sqrt{q}$ satisfying \emph{either} of the following additional 
hypotheses \cite[\emph{loc. cit.}]{Waterhouse}: (i) $(t,p) = 1$; or (ii) $p
\not \equiv 1 \mod 4$ and $t = 0$.  We will also use the fact that
the only pairs of natural numbers $(s,t)$ such that $2^s - 3^t =
\pm 1$ are $(1,0), \ (1,1) \ , \ (2,1),$ and $(3,2)$.  This follows from the 
Catalan conjecture, which has recently been proved by P. Mihailescu \cite{Mihailescu}.  We consider
several cases: \\ \\
I. If $q \in \{2, \ 4 , \ 16\}$, take $N = \ell = q+1$.  \\
II. If $q = 2^a$ for $a = 3$ or $a > 4$, take $N = q+2 = 2^a+ 2$.
If $N$ were of the form $2^b \cdot 3^c$, then $b = 1$ and
$2^{a-1}-3^c = -1$, but as recalled above there is no such $c$.
Thus $N$ is divisible by a
prime $\ell \geq 5$, and any such prime will do. \\
III. If $q = 3$, take $N = \ell = 5$. \\
IV. If $q = 3^a$ for $a > 1$, take $N = q+1$.  Since $2^b - 3^a =
1$ has no solutions for $a>1$, $N$ is divisible by a prime $\ell
\geq 5$,
and any such will do.  \\
V. If $q = p^a$ for $p \geq 5$, take $N = q-2$ and any prime $\ell
\ | \ N$.
\\ \\
Remark: When $A$ has split purely toric reduction, every class
$\eta \in H^1(K,A)$ has a unique minimal splitting extension
$L_{\eta}$, which is abelian over $K$ \cite{Gerritzen}, \cite[$\S
3.1$]{WCII}.

\section{Number fields}
\subsection{Curves of genus one}

\noindent Let $E_{/\Q}$ be the Jacobian of Selmer's cubic $3X^3 +
4Y^3 + 5Z^3 = 0$. Then $E(\Q) = 0$, $\Sha(\Q,E) \cong (\Z/3\Z)^2$
\cite{Mazur}. Note that $j(E) = 0$, so that $E$ has CM by the
maximal order in $\Q(\sqrt{-3})$.

\begin{thm}
\label{ELL}
Let $\ell$ be either $4$ or an odd prime number.  There exists a class $\eta \in H^1(\Q,E)$ such that: \\
(i) $\eta$ has period and index equal to $\ell$. \\
(ii) $\eta$ does not have an abelian splitting field.
\end{thm}
\noindent Proof: Suppose first that $\ell \geq 5$.  Then by
Poitou-Tate duality, a strong form of the local-global principle
holds in $H^1(\Q,E)[\ell^{\infty}]$, namely the natural map
\begin{equation}
\label{PT} H^1(\Q,E)[\ell^{\infty}] \ra \bigoplus_p
H^1(\Q_p,E)[\ell^{\infty}] \end{equation} is an isomorphism
\cite[I.6.26(b)]{Milne}.  (Note that $H^1(\R,E)[\ell^{\infty}] =
0$ since $\ell$ is odd.  In any case, $H^1(\R,E) = 0$ for this $E$.)  There
exist infinitely many primes $p$ such that: \\ \\
(i) $p > 3$; \\
(ii) $p \equiv -1 \pmod 3$; \\
(iii) $p \equiv -1 \pmod \ell$. \\
(iv) $E$ has good reduction mod $p$. \\ \\
Fix one such prime $p$. Condition (ii) means that $p$ is nonsplit
in the CM field $\Q(\sqrt{-3})$, so that by a well-known criterion of
Deuring, $E$ has supersingular reduction modulo $p$ (or see 
\cite[Example V.4.4]{AEC}), so that
(using (i)), $\# E(\F_p) = p+1$ and hence, by (iii), $\ell \ | \
\# E(\F_p)$.  By Lemma \ref{LOCAL2}, there exists $\eta_p \in
H^1(\Q_p,E)$ of order $\ell$.  Because the map of (\ref{PT}) is an
isomorphism, there exists a unique class $\eta \in
H^1(\Q,E)[\ell]$ restricting at $p$ to $\eta_p$ and having trivial
restriction at all other primes.  By \cite[Prop. 6]{CRELLE}, we
conclude that $\eta$ has period equals index equals $\ell$.  
Since $p \equiv -1 \pmod \ell$, evidently $\ell$ does not divide 
$p-1 = \#\mu(\Q_p)$, so by
Proposition \ref{LOCAL1}, $\eta_p$ (and \emph{a fortiori} $\eta$
itself) has no abelian splitting field.
\\ \indent
Next, note that the case $\ell = 3$ is covered by Theorem
\ref{EASY} (take $a = b = 1, \ c = 60$).  In fact it is not difficult 
to modify the above argument (taking into account that now the map of 
(\ref{PT}) has a nontrivial, but still finite, kernel) in this case, 
and one gets essentially the ``theoretical explanation'' for the 
existence of the family of curves (\ref{ONE}), since the curves 
$X^3 + pY^3 + 60p^2 Z^3$ are
indeed all principal homogeneous spaces of $E$.
\\ \indent
For $\ell = 4$, take $p = 11$ (so that $4 \ \not{|} \ p-1$); one
checks easily that $\tilde{E}(\F_{11}) \cong \Z/12\Z$.
\subsection{Curves of genus $g \geq 4$}
We will reduce to the case of curves of genus one via the
following result, a version of which was suggested to me by B.
Poonen \emph{en route} to Sabino Canyon in 2003.\footnote{This
same idea was later broached to Poonen's student S. Sharif, whose
2006 Berkeley thesis employs it as the jumping-off point for a
complete determination of the possible values of a period and
index for a genus $g$ curve over a $p$-adic field.  On the other
hand, without further assumptions on $k$ the construction of 
Proposition \ref{POONEN} is
``optimal'' in a sense that we will discuss in a later work.}
\begin{prop}
\label{POONEN} Let $K$ be a field of characteristic different from
$2$, and let $Y_{/K}$ be a genus one curve of index $n$. For any
positive integer $k$ with $kn > 1$, there exists a curve $X_{/K}$
of genus $nk+1$ and a degree two covering $X \ra Y$ defined over
$K$.
\end{prop}
\noindent Proof: By definition of the index, there exists a
$K$-irreducible divisor $D'$ on $Y$ of degree $n$; put $D_0 =
kD'$. If $n = 1$, then $D'$ consists of a single point, and since
$k > 1$, $L(|D_0|)$ is basepoint free, so there exists a divisor
$D_{\infty}$ linearly equivalent to $D_0$ and supported away from
$D_0$ and a function $f \in K(Y)$ with $\div(f) = D_0 -
D_{\infty}$.  If $n > 1$, then $L(|D'|)$ is already basepoint
free, so there exists $E'$ linearly equivalent to $D'$ and with
disjoint support, and hence a function $f \in K(Y)$ with $\div(f)
= k(D'-E')$.  The extension of function fields
$K(Y)(\sqrt{f})/K(Y)$ corresponds to a degree $2$ cover $X \ra Y$
with $2kn$ simple branch points. By the Riemann-Hurwitz theorem
$X$ has genus $kn+1$.
\\ \\
Let us now prove Theorem \ref{HIGHERG}b).  A positive integer $g
\geq 4$ may be written as $k\ell + 1$ with $k \in \Z^+$ and $\ell$
either an odd prime or $4$.  By Theorem \ref{ELL} there exists a
genus one curve $Y_{/\Q}$ of index $\ell$ without abelian points.
Applying Proposition \ref{POONEN} with $n = \ell$, we get a curve
$X$ of genus $g$ together with a degree two map $X \ra Y$. Since
$Y$ has no abelian points, neither does $X$.

\subsection{The proof of Theorem 3}
Let $n = [K:\Q]$.  Let $\ell > 7$ and $p$ be distinct primes, each
unramified in $K$, such that $\ell$ does not divide $p^a-1$ for
any $1 \leq a \leq n$; then no completion of $K$ at prime over $p$
has a rational $\ell$th root of unity. \textbf{Suppose}
$\tilde{E}_{/\F_p}$ is an elliptic curve with an $\F_p$-rational
point of order $\ell$. Lift $\tilde{E}$ to an elliptic curve
$E_{/\Q}$.  By a theorem of Ono-Skinner \cite{OS}, there exists an
elliptic curve
$E'_{/\Q}$ such that: \\ \\
(i) $E'_{/\Q_p} \cong E_{/\Q_p}$ (in particular $j(E) = j(E')$); \\
(ii) $E'$ has analytic rank zero. \\ \\
(More precisely, $E'$ is the twist of $E$ by a quadratic Dirichlet
character $\chi$ with $\chi(p) = 1$.)  By the results of
Gross-Zagier and Kolyvagin, it follows that $E'(\Q)$ and
$\Sha(\Q,E)$ are both finite.  Moreover, since $\ell > 7$, Mazur's
theorem on rational torsion points on elliptic curves \cite{Mazur}
gives $E'(\Q) \otimes \Z_{\ell} = 0$.  Now the Poitou-Tate global
duality theorem applies to show that the natural map
\[H^1(\Q,E')[\ell^{\infty}] \ra \bigoplus_p
H^1(\Q_p,E')[\ell^{\infty}] \] is a surjection.  In particular,
there exists a class $\eta \in H^1(\Q,E')[\ell^{\infty}]$ whose
local restriction $\eta_p$ has order $\ell$.  Since $K/\Q$ is
unramified at $p$, by Theorem \ref{LT}  $\eta|_{K}$ has exact
order $\ell$, and by the same arguments as above does not split
over any abelian extension of the completion of $K$ at any prime
over $p$.
\\ \indent It remains to show that, for some choices of $\ell$ and
$p$ as above, there exists $\tilde{E}_{/\F_p}$ with a point of
order $\ell$. For this, we consider primes $p > n+1$.  Another
special case of the Deuring-Waterhouse classification
\cite[Theorem 4.1]{Waterhouse} gives that for any integer $A \in
(p+1-2\sqrt{p},\ p+1+2\sqrt{p})$, there exists $E_{/\F_p}$ with $\#
E(\F_p) = A$.  In particular, for every $\ell < \sqrt{p}$, there
exists $E_{/\F_p}$ with an element of order $\ell$.  There are at
most $n$ primes $\ell < \sqrt{p}$ whose order, as elements of
$\F_p^{\times}$, is at most $n$, so when $p$ is sufficiently large
compared to $n$ there are many primes $\ell$ such that there
exists $E_{/\F_p}$ with elements of order $\ell$. This completes
the proof of the theorem. \noindent
\\ \\
Remark: The proof we have given is a veritable showcase of the
deepest results of \emph{fin de si\`ecle} elliptic curve theory.
More care would probably lead to a more elementary (but perhaps
less amusing) proof.
\\ \\
%

\section{Some final remarks}
\noindent With a single exception, this section contains not
results but connections to other work, advertisements,
conjectures, questions, and even ``hearsay.''

\subsection{Conjectural strengthenings}
I find it likely that the following stronger statements hold:
\begin{conj}
Let $K$ be any number field and $E_{/K}$ any elliptic curve.\\
a) There is a genus one curve
$C_{/K}$, with Jacobian $E$, and $C(K^{\ab}) = \emptyset$. \\
b) For all $d \geq 3$, there is a degree $d$ plane curve
$C_{/\Q}$ such that $C(K^{\ab}) = \emptyset$. \\
c) For all $g \geq 4$, there is a curve $C_{/\Q}$ of genus $g$
such that $C(K^{\ab}) = \emptyset$. \\
\end{conj}
\noindent

\subsection{Higher-dimensional varieties}
The following immediate generalization of Theorem \ref{EASY}
gives, for any odd prime $\ell$, a Calabi-Yau $(\ell-2)$-fold
$V_{/\Q}$ without $\Q^{\ab}$-rational points.
\begin{prop}
\label{CY} Let $\ell$ be an odd prime and $p$ a prime with $\ell
\not{|} \ p(p-1)$.  Then
\[\sum_{i=0}^{\ell-1} p^i X_i^{\ell} = 0 \]
has no points over $\Q_p^{\ab}$.
\end{prop}
\noindent The proof does not go through for composite $\ell$: in
particular for $\ell = 4$ the construction does not yield examples
of K3 surfaces without abelian points. Note that, given a quartic
surface without abelian points, taking a general hyperplane
section would give a negative answer to Question \ref{GENUS3}.
\\ \\
It would be of great interest to find a geometrically rational
variety $V_{/\Q}$ without abelian points.  Especially, it is a
famous open question of Artin whether a Fano hypersurface (i.e.,
the zero locus of a degree $d$ homogeneous polynomial in at least
$d+1$ variables) defined over $\Q^{\ab}$ must have a rational
point.  As mentioned above, by a theorem of Lang \cite{Lang}, such
hypersurfaces have points everywhere locally.  By work of Brauer
and Birch, the existence of quadratic points is known for
``sufficiently Fano'' hypersurfaces: for every fixed $d$, there
exists $n = n(d)$ such that a degree $d$ form in $n$ variables
over $\Q$ has a nontrivial solution in (e.g.) $\Q(\sqrt{-1})$.
Finally, work of Kanevsky \cite{Kanevsky} shows that if $V_{/K}$ is a cubic 
surface over a number field $K$ such that for all finite extensions 
$L/K$, the Brauer-Manin obstruction to the existence of 
$L$-rational points on $V_{/L}$ is the only one, then $V(K^{\ab}) 
\neq \emptyset$.

\noindent \subsection{Varieties \emph{with} abelian points} One
can ask for nontrivial examples of varieties with abelian points,
i.e., not coming from a quadratic covering of a variety with
$\Q$-points. The best example I know is that of Severi-Brauer
varieties (varieties $V_{/\Q}$ such that $V_{/\overline{\Q}} \cong
\PP^N_{/\overline{\Q}}$): there are always abelian points, but, in 
dimension at least two, usually not quadratic points.  
Indeed, the
Brauer-Hasse-Noether theorem says that every element of the Brauer
group of a number field is given by a cyclic algebra.\footnote{Or
see \cite[$\S$ II.3.3, Prop. 9]{CG} for a more elementary argument
showing that $\Br(\Q^{\ab}) = 0$.}
\\ \\
Somewhat distressingly, the following seems to be nontrivial:
\begin{prob}
\label{PROB} For all positive integers $n$ (or even for infinitely
many $n$), exhibit a genus one curve $C_{/\Q}$ of index $n$ and
such that $C(\Q^{\ab}) \neq \emptyset$.
\end{prob}
\noindent Perhaps the solution to Problem \ref{PROB} will involve
some Iwasawa theory.
\subsection{Varieties with points everywhere locally}
\noindent A motivation for this work came from a question of D.
Jetchev, who asked whether Selmer's curve $3X^3 + 4Y^3 + 5Z^3 = 0$
has points in any abelian cubic field, or in any abelian number
field.  Theorem \ref{EASY} gives negative answers for
superficially similar curves: e.g. the curve $2X^3 + 4Y^3 + 5Z^3 =
0$ fails to have abelian points.  But Selmer's curve has points
everywhere locally, rendering useless the present approach.
\\ \\
Rephrasing the question slightly, we may ask:
\begin{ques}
\label{JETCHEV} Fix a positive integer $n$.  Does every locally
trivial genus one curve of index $n$ defined over a number field
have an abelian point?
\end{ques}
\noindent I have been told that random matrix theory predicts a
positive answer when $n = 3$.

\subsection{Solvable points}
We were also motivated by work in progress of M. Ciperiani and A.
Wiles, who study \emph{solvable} points on curves of genus one.
They are able to show (at least) that a genus one curve $C_{/\Q}$
which is locally trivial and with semistable Jacobian has a
solvable point.
\\ \\
Using the solvability of the absolute Galois groups of $\Q_p$ (and
$\R$), it is easy to see that for every variety $V$ over a number
field $K$, there is a solvable extension $L/K$ such that $V_{/L}$
has points everywhere locally.  Thus, an affirmative answer to
Question \ref{JETCHEV} for all $n$ (which I must admit seems 
unlikely) would imply the existence of
solvable points on all curves of genus one. 
\subsection{Metabelian points}
Note that in Section 2, all our examples of principal homogeneous
spaces over $\Q_p$ without abelian points have points over the
maximal abelian extension of $\Q_p^{\ab}$, i.e., over a
\textbf{metabelian} extension of $\Q_p$.  It was suggested to me a
few years ago by B. Mazur that every genus one curve over $\Q$
should have metabelian points.  As far as I know, this remains
open even over $\Qp$, although special cases follow from results
of Lang-Tate \cite{LT} and Lichtenbaum \cite{Lichtenbaum}.

\section*{Appendix: varieties over $\Q((t))$} A somewhat different perspective would
be to fix a ``class'' of algebraic varieties (e.g., curves of a
given genus $g$, or hypersurfaces of degree $d$ in $\PP^N$) and
ask whether for \emph{any} field $F$,\footnote{Just for
simplicity, let us assume that $F$ has characteristic $0$.} a
variety $V_{/F}$ of this type must have points in the maximal
abelian extension of $F$. With ``abelian'' replaced by
``solvable,'' this is the setting of recent work of A. P\'al
\cite{Pal}. In this appendix we will show that, with a suitable
choice of $F$, there are additional classes of $F$-varieties
without $F^{\ab}$-points.
\\ \indent
%
\newcommand{\abc}{\operatorname{abc}}
\newcommand{\tonehalf}{t^{\frac{1}{2}}}
\noindent It is convenient to work with $F = \Q((t))$.  The
absolute Galois group $\mathfrak{g}_F$ of $F$ lies in a split
exact sequence
\[1 \ra \hat{\Z} \ra \mathfrak{g}_F \ra \mathfrak{g}_{\Q} \ra 1 \]
where the action of $\mathfrak{g}_{\Q}$ on $\hat{\Z}$ is by the
cyclotomic character.  It follows that the maximal abelian
extension of $F$ is generated by the roots of unity together with
$\tonehalf$.  The field $F(\mu_{\infty},t^{\frac{1}{2}})$ is
Henselian with respect to the discrete valuation
$\frac{\ord_t}{2}$.  We will work instead with its completion
$\Q^{\ab}((t^{\frac{1}{2}}))$, which -- by a small abuse of
notation -- we will denote by $F^{\ab}$.  Since this field contains 
what is literally the maximal abelian extension of $F$, finding 
varieties $V_{/F}$ with $V(F^{\ab}) = \emptyset$ is \emph{a priori} 
a stronger result than showing that they do not have abelian points.  
(But in fact it is equivalent: a variety
defined over a discretely valued field has points rational over the
Henselization iff it has points rational over the completion.)

%
\begin{prop}
Every Severi-Brauer variety $V_{/F}$ has an abelian point, but
there exist Severi-Brauer varieties $V_{/F^{\ab}}$ without
rational points.
\end{prop}
\noindent Proof: This may be viewed as a question about the Brauer
groups $\Br(F)$ and $\Br(F^{\ab})$ (e.g. \cite[$\S X.6$]{CL}). For
any complete, discretely valued field $K$ with perfect residue
field $K$, there is an exact sequence
\[0 \ra \Br K \ra \Br K \stackrel{c}{\ra} X(\mathfrak{g}_k) \ra 0, \]
where the last term is the character group of the Galois group of
the residue field \cite[Theorem X.3.2]{CL}.  Consider first the
case of $K = \Q((t))$, $k = \Q$.  Then, for $\alpha \in \Br K$ we
can split the character $c(\alpha)$ via a unique unramified
abelian extension $L/K$ with (abelian) residue extension $l/\Q$.
By the exact sequence, $\alpha|_{L} \in \Br(l)$.  Now every
element of the Brauer group of a number field can be split by a
cyclotomic extension, so overall we get that $\Br(\Q((t))) =
\Br(\Q^{\ab}((t))/\Q((t)))$, giving the first statement of the
proposition. \\ \indent On the other hand, $F^{\ab} =
\Q^{\ab}((t^{1/2})) \cong \Q^{\ab}((t))$ is again a local field,
whose residue field $\Q^{\ab}$ has trivial Brauer group but highly
nontrivial character group.\footnote{We are again encountering the
phenomenon that $(\Q^{\ab})^{\ab}$, the maximal metabelian
extension of $\Q$, is very much larger than $\Q^{\ab}$.}  For
example, the quaternion algebra $\langle \sqrt{2}, \tonehalf
\rangle$ represents a nontrivial element of
$\Br(F^{\ab})$. %
\begin{cor}
There exists a geometrically rational $4$-fold $V_{/F}$ with
$V(F^{\ab}) = \emptyset$.
\end{cor}
\noindent Proof: The quaternion algebra $\langle \sqrt{2}, \
\tonehalf \rangle$ defined over $F' = F(\sqrt{2},\tonehalf)$ is
nonsplit over $F^{\ab}$; this corresponds to a conic $C_{/F'}$
with $C(F^{\ab}) = \emptyset$.  Restriction of scalars (or ``Weil
restriction'') from $F'$ to $F$ gives a fourfold $V_{/F}$ such
that $V_{/\overline{F}} \cong (\PP^1)^4$ and $V(F^{\ab}) = C(F'
\otimes_F F^{\ab}) = C(\prod_{i=1}^4 F^{\ab}) = C(F^{\ab})^4 =
\emptyset$.
\\ \\
Remark: A famous theorem of Merkurjev implies that over an
arbitrary field $K$, Severi-Brauer varieties $V_{/K}$ have
metabelian points.  Work of Wedderburn and Albert shows that every
Severi-Brauer surface is split over a cyclic cubic extension and
every Severi-Brauer threefold is split over a $\Z/2\Z
\times\Z/2\Z$-extension.  I am not personally in possession of an
example of a division algebra without an abelian splitting field,
but I presume this is the generic situation when the index is
divisible by a sufficiently high power of a prime.
\\ \\
Recall the \textbf{norm form} $N(L/K)$ associated to a separable
field extension of degree $d$: choosing a $K$-basis
$\alpha_1,\ldots,\alpha_d$ of $l$, the map \[ N: K^d \ra K, \ 
(x_1,\ldots,x_d) \mapsto N(\sum_{i=1}^d x_i \alpha_i) \] is a
degree $d$ polynomial. Let us say that a homogeneous $K$-form
$f(X_1,\ldots,X_d)$ is \emph{isotropic} (resp. \emph{anisotropic})
if there exist $(x_1,\ldots,x_n)$, not all zero, with
$f(x_1,\ldots,x_d) = 0$ (resp. there is only the zero solution).
\\ \\
If $M/K$ is a separable field extension, then $N(L/K)_{/M}$ is the
norm form for the extension of algebras $L \otimes_K M / M$, and is 
anisotropic if and only if $L \otimes_K M$ is a field. In other
words:
\begin{lemma}
\label{NORM} Let $M/K$ be a separable field extension.  Then the
norm form $N = N(L/K)$ is anistropic over $M$ if and only if $M$
and $L$ are linearly disjoint over $K$.  In particular, if $L/K$
is not abelian, then $N(L/K)$ is anisotropic over every abelian
extension of $K$.
\end{lemma}
\noindent Note however that the norm form $N(L/K)$ of a nontrivial
field extension is geometrically reducible.  Indeed, over the
Galois closure of $L/K$, with a suitable choice of basis $N$
becomes $N(X_1,\ldots,X_d) = X_1\cdots X_d$ ($d = [L:K]$).  The
corresponding closed subscheme is of dimension $d-2$ and has as
its singular locus a finite union of linear subspaces of dimension
$d-3$.
\begin{prop}
Fix a positive integer $d \geq 3$, and let $K/\Q$ be a degree $d$
number field with Galois group $S_d$.  Let $N$ be the norm form of
$K/\Q$.  Then
\begin{equation}
\label{NORM} N(X_1,\ldots,X_d) = t Z^n \end{equation} has no
$F^{\ab}$-rational points.
\end{prop}
\noindent Proof: There are no solutions with $Z = 0$ since $K(T) =
K \otimes_{\Q} F / F$ is linearly disjoint from $F^{\ab}$. A
solution with $Z \neq 0$ exists if and only if $t$ is a norm from
$M = K \otimes F^{\ab}$ down to $L = F^{\ab}$.  But $M/L$ is a
finite unramified extension of complete discretely valued fields,
so the image of the norm map consists of elements whose valuation
is divisible by $d$, whereas $v_L(t) = 2$.
\begin{cor}
There is a geometrically rational surface $S_{/F}$ without abelian
points.
\end{cor}
\noindent Proof: Taking $d = 3$, one gets a geometrically integral
cubic suface $S_{/F}$ without abelian points, but with finitely
many (in fact $3$) singular points.  We can resolve the
singularities by a birational $F$-morphism $\tilde{S} \ra S$, and
$S(F^{\ab}) = \emptyset$ implies $\tilde{S}(F^{\ab}) = \emptyset$.
\begin{cor}
There is a plane curve of any degree $d \geq 3$ without abelian
points.
\end{cor}
\noindent Proof: Since the singular locus of (\ref{NORM}) is a
finite union of codimension two affine subspaces, one sees easily
that intersecting with a general $2$-plane gives a nonsingular
curve. In particular, taking $n = 4$ we get a genus $3$ curve
$C_{/F}$ with $C(F^{\ab}) = \emptyset$.
\\ \\
Following \cite{Pal}, we get another approach to curves of genus
$3$ without abelian points:
\begin{prop}
Let $F$ be a complete, discretely valued field whose residue field
$f$ contains an extension $m$ which is Galois with group
isomorphic to $S_4$.  Then there exists a genus $3$ curve $C_{/F}$
with $C(F^{\ab}) = \emptyset$.
\end{prop}
\noindent Proof: After choosing an isomorphism of $G = \Gal(m/f)$
with $S_4$, we get an action of $\mathfrak{g}_f$ on the complete
graph $K_4$ on $4$ vertices.  By \cite[Prop. 4.6]{Pal}, there
exists a stable curve $C_{/f}$ with rational
geometric components, whose corresponding dual graph is
isomorphic, as a $\mathfrak{g}_{f}$-module, to $K_4$ with the
chosen $G$-action. By \cite[Cor. 4.4]{Pal}, there exists a stable
curve over the valuation ring $R_F$ of $F$ whose generic fiber is an
honest (i.e., nonsingular and geometrically integral) curve
$C_{/F}$ and whose special fiber is isomorphic to $C_{/f}$.
Since the stabilizer of any vertex or edge of $K_4$ is a non-normal 
subgroup of $S_4$, after making any abelian residue extension $f'/f$, there 
are no $\mathfrak{g}_{f'}$-fixed vertices or edges of the dual graph, 
so $C(f') = \emptyset$.  Hence if $F'/F$ is any extension with 
abelian residue extension -- so in particular if $F'/F$ is itself 
abelian -- $C(f') = \emptyset$ implies $C(R_{F'}) = C(F') = \emptyset$.


\begin{thebibliography}{X}
%
\bibitem[1]{CRELLE} P.L. Clark, \emph{There are genus one curves
of every index over every number field}, to appear in J. Reine
Angew. Math.
%
\bibitem[2]{WCII} P.L. Clark, \emph{Period-index problems in
WC-groups II: abelian varieties}, submitted.
%
\bibitem[3]{FJ} M. Fried and M. Jarden, \emph{Field Arithmetic},
Ergebnisse der Mathematik (3) 11 (1986), Springer-Verlag.
%
\bibitem[4]{Frey} G. Frey, \emph{Pseudo algebraically closed
fields with non-Archimedean real valuations}, J. Algebra 26
(1973), 202-207.
%
\bibitem[5]{Gerritzen} L. Gerritzen.
\emph{Periode und Index eines prinzipal-homogenen Raumes \"uber
gewissen abelschen Variet\"aten}, Manuscripta Math. 8 (1973),
131-142.
%
\bibitem[6]{Kanevsky} D. Kanevsky, \emph{Application of the 
conjecture on the Manin obstruction to various Diophantine 
problems}, Ast\'erisque 147-148 (1987), 307-314, 345.
%
\bibitem[7]{Lang} S. Lang, \emph{On quasi-algebraic closure}, Ann.
of Math. 55 (1952), 373-390.
%
\bibitem[8]{LT} S. Lang and J. Tate, \emph{Principal homogeneous
spaces over abelian varieties}, Amer. J. Math. 80 (1958), 659-684.
%
\bibitem[9]{Lichtenbaum} S. Lichtenbaum, \emph{The period-index
problem for elliptic curves}, Amer. J. Math. 90 (1968), 1209-1223.
%

\bibitem[10]{Mazur} B. Mazur, \emph{On the passage from local to
global in number theory}, Bull. Amer. Math. Soc. (N.S.) 29 (1993),
14-50.


\bibitem[11]{Mihailescu} P. Mihailescu, \emph{Primary cyclotomic units and
a proof of Catalan's conjecture}, J. Reine Angew. Math. 572
(2004), 167-195.
%
\bibitem[12]{Milne} J. Milne.
\emph{Arithmetic Duality Theorems}, Perspectives in Mathematics,
1.  Academic Press Inc., 1986.
%
%
\bibitem[13]{OS} K. Ono and C. Skinner, \emph{Non-vanishing of
quadratic twists of modular $L$-functions}, Invent. Math. 134
(1998), 651-660.
%
\bibitem[14]{Pal} A. P\'al, \emph{Solvable points on projective
algebraic curves}, Canad. J. Math. 56 (2004), 612-637. 
%
\bibitem[15]{CL} J.-P. Serre, \emph{Corps locaux}, Hermann, Paris,
1962.
%
\bibitem[16]{CG} J.-P. Serre, \emph{Galois Cohomology}, Lecture
Notes in Mathematics 5, 5th revised edition, Springer-Verlag,
1994.
%
\bibitem[17]{AEC} J. Silverman, \emph{The arithmetic of elliptic 
curves}, Graduate Texts in Mathematics 106, Springer-Verlag, 1986.
%
\bibitem[18]{Tate} J. Tate, \emph{WC-groups over $p$-adic fields},
Sem. Bourbaki, Exp. 156, 1957.
%
\bibitem[19]{Waterhouse} W.C. Waterhouse, \emph{Abelian varieties
over finite fields}, Ann. Sci. \`Ecole Norm. Sup. 2 (1969), 521-560.


\end{thebibliography}
\end{document}